\documentclass[final,onefignum,onetabnum]{siamart171218}

\usepackage{lipsum}
\usepackage{amsfonts}
\usepackage{graphicx}
\usepackage{epstopdf}
\usepackage{algorithmic}
\usepackage{url,tikz,verbatim,booktabs,amsmath,bm,amssymb,subfigure}

\ifpdf
  \DeclareGraphicsExtensions{.eps,.pdf,.png,.jpg}
\else
  \DeclareGraphicsExtensions{.eps}
\fi

\newcommand{\mat}[1]{\bm{#1}}
\newcommand{\ten}[1]{\bm{\mathcal{#1}}}

\newsiamremark{remark}{Remark}
\newsiamremark{hypothesis}{Hypothesis}
\crefname{hypothesis}{Hypothesis}{Hypotheses}
\newsiamthm{claim}{Claim}

\headers{Tensor Ring Trouble}{K. Batselier}

\title{The trouble with tensor ring decompositions}

\author{Kim Batselier\thanks{Delft Center for Systems and Control, Delft University of Technology, Delft, The Netherlands
  (\email{k.batselier@tudelft.nl}, \url{https://sites.google.com/view/kim-batselier/home}).}}

\usepackage{amsopn}

\ifpdf
\hypersetup{
  pdftitle={Tensor Ring Trouble},
  pdfauthor={K. Batselier}
}
\fi

\begin{document}

\maketitle

\begin{abstract}
The tensor train decomposition decomposes a tensor into a ``train" of 3-way tensors that are interconnected through the summation of auxiliary indices. The decomposition is stable, has a well-defined notion of rank and enables the user to perform various linear algebra operations on vectors and matrices of exponential size in a computationally efficient manner. The tensor ring decomposition replaces the train by a ring through the introduction of one additional auxiliary variable. This article discusses a major issue with the tensor ring decomposition: its inability to compute an exact minimal-rank decomposition from a decomposition with sub-optimal ranks. Both the contraction operation and Hadamard product are motivated from applications and it is shown through simple examples how the tensor ring-rounding procedure fails to retrieve minimal-rank decompositions with these operations. These observations, together with the already known issue of not being able to find a best low-rank tensor ring approximation to a given tensor indicate that the applicability of tensor rings is severely limited.
\end{abstract}

\begin{keywords}
  tensors, tensor ring, tensor decomposition, minimal-rank
\end{keywords}

\begin{AMS}
 15A23, 15A69, 65F99
\end{AMS}

\section{Introduction}
The recent development of the tensor train (TT) decomposition~\cite{Oseledets2010,ivanTT,ttcross} in the domain of scientific computing led to myriad applications. These tensor decompositions are a tool to lift the infamous curse of dimensionality where the key idea is to represent a given tensor (or vector of exponential length) by its TT or TR-decomposition. The hope is then that a low-rank representation results in either an exact representation or in an approximation that is sufficient for practical purposes. The innate power and versatility of this decomposition to lift the curse of dimensionality has led to applications in signal processing~\cite{MAL-067}, neural networks~\cite{LGROL14,Novikov2015}, supervised machine learning~\cite{TNClassifier,stoudenmire2016supervised,stoudenmire2018learning}, and many others~\cite{MAL-059,Lee2015,Lee2016,Os-mvk2-2011,Oseledets2012}. The TT-decomposition was in fact discovered a decade earlier in the physics community, where it was known as the Matrix Product State (MPS) ansatz~\cite{MPS1997,SCHOLLWOCK201196}. There are two versions of the MPS: one with open boundary conditions and one with periodic boundary conditions. The open boundary conditions case is what is now known as a TT-decomposition, while the periodic boundary conditions case is commonly referred to as the Tensor Ring (TR) decomposition. The generalization of TT to TR happened very shortly after the introduction of TTs into the scientific computing community~\cite{Espig2011,Espig2012,Khoromskij2011}. In spite of their generality, TRs have found only limited application: in the compression of fully connected and convolutional layers in neural networks~\cite{Wang_TRcompletion} and in solving the tensor completion problem~\cite{wang2017efficient,yuan2018rank}.

From an engineering point of view, the compression of a large data set (e.g. by means of a TT or TR) is usually not the end goal but rather a first step in the design and implementation of a solution to a particular problem. The next step is then usually to perform particular operations on the compressed data in order to compute some desired result. This in part explains the power and versatility of the TT/TR representation in that they allows us ``to do linear algebra" with vectors and matrices of exponential size. Indeed, both vectors and matrices can be represented in either an exact or approximate manner by TTs/TRs. Linear algebra operations such as addition, matrix  multiplication,  the Hadamard and inner product, the Fast Fourier Transform~\cite{dolgov2012superfast}, linear system solving and matrix inversion~\cite{Oseledets2012} and the computation of low-rank SVD approximations~\cite{TNrSVD,Lee2015,Lee2016} can then be performed on the TT/TR representation directly, which often results in enormous savings in computation and storage. There is a small caveat, however, in that many of these operations will result in TTs/TRs with increased sub-optimal ranks. This issue is completely solved by the rounding procedure for TT-decompositions~\cite[p.~2305]{ivanTT}. The TT-rounding algorithm reduces any sub-optimal TT-ranks to their minimal values through subsequent QR and SVD factorizations and has been recently adapted to work for TRs~\cite{mickelin2018tensor}. This article will show how the TR-rounding algorithm fails to recover the minimal TR-ranks for various cases. This issue is different from the already known problem that a best low-rank TR approximation of a given tensor may not exist~\cite{Landsburg2012,ye2018tensor}, as only exact minimal-rank decompositions are considered in this article. Before illustrating the problem in Section~\ref{sec:minrank} we first provide a short summary on the TR decomposition in Section~\ref{sec:TR}.

\section{Tensor ring decomposition}
\label{sec:TR}
A $d$-way tensor $\ten{A} \in\mathbb{R}^{I_1\times I_2\times \cdots\times I_d}$ is in this article a $d$-dimensional array where each entry is completely determined by $d$ indices $i_1,i_2,\ldots,i_d$. The convention $1\leq i_k\leq I_k$ for $k=1,\ldots,d$ is used, together with MATLAB notation to denote entries of tensors. Boldface capital calligraphic letters $\ten{A},\ten{B},\ldots$ are used to denote tensors, boldface capital letters $\mat{A},\mat{B},\ldots$ denote matrices, boldface letters $\mat{a},\mat{b},\ldots$ denote vectors, and Roman letters $a,b,\ldots$ denote scalars. The Hadamard (elementwise) product of two tensors is denoted $\circ$. The definition of the TR decomposition has appeared in~\cite{Espig2011,Espig2012,Khoromskij2011,mickelin2018tensor,zhao2016tensor} and is presented here for completeness.
\begin{definition}
\label{def:TR}
The TR decomposition of a given tensor  $\ten{A} \in\mathbb{R}^{I_1\times I_2\times \cdots\times I_d}$ is a set of 3-way tensors $\ten{A}^{(k)} \in \mathbb{R}^{R_k \times I_k \times R_{k+1}}\; (1 \leq k \leq d)$ such that each entry $\ten{A}(i_1,i_2,\ldots,i_d)$ can be computed from
\begin{align}
     & \sum_{r_1=1}^{R_1} \sum_{r_2=1}^{R_2} \cdots \sum_{r_d=1}^{R_d} \ten{A}^{(1)}(r_1,i_1,r_2)\,\ten{A}^{(2)}(r_2,i_2,r_3)\, \ldots, \,\ten{A}^{(d)}(r_d,i_d,r_1).
     \label{eq:TRdef}
\end{align}
\end{definition}
The summations in \eqref{eq:TRdef} need to return a scalar, which implies that $R_{d+1}=R_1$. The minimal dimensions $R_1,R_2,\ldots,R_d$ of the auxiliary indices $r_1,r_2,\ldots,r_d$ such that each tensor entry $\ten{A}(i_1,i_2,\ldots,i_d)$ is equal to \eqref{eq:TRdef} are per definition the TR-ranks. Fixing the $i_1,i_2,\ldots,i_d$ indices to definite values turns each of the TR-tensors into a matrix $\mat{A}^{(k)} \triangleq \ten{A}^{(k)}(:,i_k,:) \in \mathbb{R}^{R_k \times R_{k+1}}$, allowing us to rewrite~\eqref{eq:TRdef} as
\begin{align}
\ten{A}(i_1,i_2,\ldots,i_d) &= \textrm{Trace} \left( \mat{A}^{(1)} \mat{A}^{(2)} \cdots \mat{A}^{(d)}\right),
\label{eq:TRdefalt}
\end{align}
where the summations over the auxiliary indices in~\eqref{eq:TRdef} are now interpreted as matrix matrix multiplications together with the trace operation. The TT-decomposition is commonly defined as a TR for which $R_1=1$. Strictly speaking, the location of the unit-rank in the ring is of no significance as the cyclic permutation property of the trace operator in~\eqref{eq:TRdefalt} always allows us to move the unit-rank to both the first and last matrix factor. Consequently, an alternative definition of the TT-decomposition is then Definition~\ref{eq:TRdef} with the additional condition that at least one TR-rank is one. Any TR can be converted into a TT by rewriting the TR as a sum of $r_k$ TTs, where $r_k$ is the auxiliary index that is removed~\cite[p.~21]{mickelin2018tensor}. An alternative method to convert a TR into a TT is presented in~\cite[p.~2]{handschuh2012changing}. An algorithm for the conversion of a TT into a TR is described in~\cite[p.22]{mickelin2018tensor}.

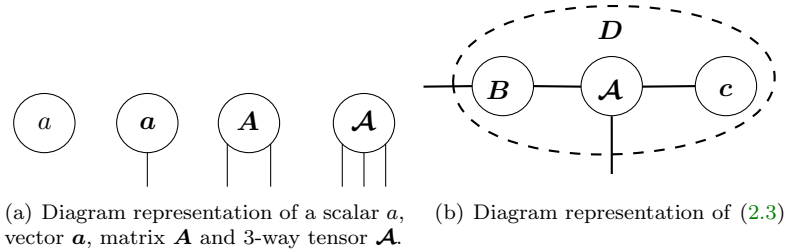
\begin{figure}[h]
\centering
\subfigure[Diagram representation of a scalar $a$, vector $\mat{a}$, matrix $\mat{A}$ and 3-way tensor $\ten{A}$.]{
\begin{minipage}[b]{0.4\linewidth}
\ifx\du\undefined
  \newlength{\du}
\fi
\setlength{\du}{4\unitlength}
\begin{tikzpicture}
\pgftransformxscale{1.000000}
\pgftransformyscale{-1.000000}
\definecolor{dialinecolor}{rgb}{0.000000, 0.000000, 0.000000}
\pgfsetstrokecolor{dialinecolor}
\definecolor{dialinecolor}{rgb}{1.000000, 1.000000, 1.000000}
\pgfsetfillcolor{dialinecolor}
\definecolor{dialinecolor}{rgb}{1.000000, 1.000000, 1.000000}
\pgfsetfillcolor{dialinecolor}
\pgfpathellipse{\pgfpoint{-6.073446\du}{10.779091\du}}{\pgfpoint{2.900000\du}{0\du}}{\pgfpoint{0\du}{2.800000\du}}f
\pgfusepath{fill}
\pgfsetlinewidth{0.100000\du}
\pgfsetdash{}{0pt}
\pgfsetdash{}{0pt}
\definecolor{dialinecolor}{rgb}{0.000000, 0.000000, 0.000000}
\pgfsetstrokecolor{dialinecolor}
\pgfpathellipse{\pgfpoint{-6.073446\du}{10.779091\du}}{\pgfpoint{2.900000\du}{0\du}}{\pgfpoint{0\du}{2.800000\du}}
\pgfusepath{stroke}
\pgfsetlinewidth{0.100000\du}
\pgfsetdash{}{0pt}
\pgfsetdash{}{0pt}
\pgfsetbuttcap
{
\definecolor{dialinecolor}{rgb}{0.000000, 0.000000, 0.000000}
\pgfsetfillcolor{dialinecolor}
\definecolor{dialinecolor}{rgb}{0.000000, 0.000000, 0.000000}
\pgfsetstrokecolor{dialinecolor}
\draw (-4.022836\du,12.758990\du)--(-4.022836\du,16.80\du);
}
\pgfsetlinewidth{0.100000\du}
\pgfsetdash{}{0pt}
\pgfsetdash{}{0pt}
\pgfsetbuttcap
{
\definecolor{dialinecolor}{rgb}{0.000000, 0.000000, 0.000000}
\pgfsetfillcolor{dialinecolor}
\definecolor{dialinecolor}{rgb}{0.000000, 0.000000, 0.000000}
\pgfsetstrokecolor{dialinecolor}
\draw (-8.124055\du,12.758990\du)--(-8.124055\du,16.80\du);
}
\pgfsetlinewidth{0.100000\du}
\pgfsetdash{}{0pt}
\pgfsetdash{}{0pt}
\pgfsetbuttcap
{
\definecolor{dialinecolor}{rgb}{0.000000, 0.000000, 0.000000}
\pgfsetfillcolor{dialinecolor}
\definecolor{dialinecolor}{rgb}{0.000000, 0.000000, 0.000000}
\pgfsetstrokecolor{dialinecolor}
\draw (-6.073446\du,13.579091\du)--(-6.073446\du,16.80\du);
}
\definecolor{dialinecolor}{rgb}{0.000000, 0.000000, 0.000000}
\pgfsetstrokecolor{dialinecolor}
\node[anchor=west] at (-8.20\du,10.6\du){$\ten{A}$};
\definecolor{dialinecolor}{rgb}{1.000000, 1.000000, 1.000000}
\pgfsetfillcolor{dialinecolor}
\pgfpathellipse{\pgfpoint{-36.297515\du}{10.779091\du}}{\pgfpoint{2.900000\du}{0\du}}{\pgfpoint{0\du}{2.800000\du}}
\pgfusepath{fill}
\pgfsetlinewidth{0.100000\du}
\pgfsetdash{}{0pt}
\pgfsetdash{}{0pt}
\definecolor{dialinecolor}{rgb}{0.000000, 0.000000, 0.000000}
\pgfsetstrokecolor{dialinecolor}
\pgfpathellipse{\pgfpoint{-36.297515\du}{10.779091\du}}{\pgfpoint{2.900000\du}{0\du}}{\pgfpoint{0\du}{2.800000\du}}
\pgfusepath{stroke}
\definecolor{dialinecolor}{rgb}{0.000000, 0.000000, 0.000000}
\pgfsetstrokecolor{dialinecolor}
\node[anchor=west] at (-37.80\du,10.80\du){$a$};
\definecolor{dialinecolor}{rgb}{1.000000, 1.000000, 1.000000}
\pgfsetfillcolor{dialinecolor}
\pgfpathellipse{\pgfpoint{-26.560105\du}{10.779091\du}}{\pgfpoint{2.900000\du}{0\du}}{\pgfpoint{0\du}{2.800000\du}}
\pgfusepath{fill}
\pgfsetlinewidth{0.100000\du}
\pgfsetdash{}{0pt}
\pgfsetdash{}{0pt}
\definecolor{dialinecolor}{rgb}{0.000000, 0.000000, 0.000000}
\pgfsetstrokecolor{dialinecolor}
\pgfpathellipse{\pgfpoint{-26.560105\du}{10.779091\du}}{\pgfpoint{2.900000\du}{0\du}}{\pgfpoint{0\du}{2.800000\du}}
\pgfusepath{stroke}
\pgfsetlinewidth{0.100000\du}
\pgfsetdash{}{0pt}
\pgfsetdash{}{0pt}
\pgfsetbuttcap
{
\definecolor{dialinecolor}{rgb}{0.000000, 0.000000, 0.000000}
\pgfsetfillcolor{dialinecolor}
\definecolor{dialinecolor}{rgb}{0.000000, 0.000000, 0.000000}
\pgfsetstrokecolor{dialinecolor}
\draw (-26.560105\du,13.579091\du)--(-26.560105\du,16.80\du);
}
\definecolor{dialinecolor}{rgb}{0.000000, 0.000000, 0.000000}
\pgfsetstrokecolor{dialinecolor}
\node[anchor=west] at (-28.20\du,10.80\du){$\mat{a}$};
\definecolor{dialinecolor}{rgb}{1.000000, 1.000000, 1.000000}
\pgfsetfillcolor{dialinecolor}
\pgfpathellipse{\pgfpoint{-16.948419\du}{10.779091\du}}{\pgfpoint{2.900000\du}{0\du}}{\pgfpoint{0\du}{2.800000\du}}
\pgfusepath{fill}
\pgfsetlinewidth{0.100000\du}
\pgfsetdash{}{0pt}
\pgfsetdash{}{0pt}
\definecolor{dialinecolor}{rgb}{0.000000, 0.000000, 0.000000}
\pgfsetstrokecolor{dialinecolor}
\pgfpathellipse{\pgfpoint{-16.948419\du}{10.779091\du}}{\pgfpoint{2.900000\du}{0\du}}{\pgfpoint{0\du}{2.800000\du}}
\pgfusepath{stroke}
\pgfsetlinewidth{0.100000\du}
\pgfsetdash{}{0pt}
\pgfsetdash{}{0pt}
\pgfsetbuttcap
{
\definecolor{dialinecolor}{rgb}{0.000000, 0.000000, 0.000000}
\pgfsetfillcolor{dialinecolor}
\definecolor{dialinecolor}{rgb}{0.000000, 0.000000, 0.000000}
\pgfsetstrokecolor{dialinecolor}
\draw (-14.897810\du,12.758990\du)--(-14.897810\du,16.80\du);
}
\pgfsetlinewidth{0.100000\du}
\pgfsetdash{}{0pt}
\pgfsetdash{}{0pt}
\pgfsetbuttcap
{
\definecolor{dialinecolor}{rgb}{0.000000, 0.000000, 0.000000}
\pgfsetfillcolor{dialinecolor}
\definecolor{dialinecolor}{rgb}{0.000000, 0.000000, 0.000000}
\pgfsetstrokecolor{dialinecolor}
\draw (-18.999029\du,12.758990\du)--(-18.999029\du,16.80\du);
}
\definecolor{dialinecolor}{rgb}{0.000000, 0.000000, 0.000000}
\pgfsetstrokecolor{dialinecolor}
\node[anchor=west] at (-19.00\du,10.6\du){$\mat{A}$};
\end{tikzpicture}
\end{minipage}}
\subfigure[Diagram representation of~\eqref{eqn:simpleTN} ]{
\begin{minipage}[b]{0.4\linewidth}
\ifx\du\undefined
  \newlength{\du}
\fi
\setlength{\du}{4\unitlength}
\begin{tikzpicture}
\pgftransformxscale{1.000000}
\pgftransformyscale{-1.000000}
\definecolor{dialinecolor}{rgb}{0.000000, 0.000000, 0.000000}
\pgfsetstrokecolor{dialinecolor}
\definecolor{dialinecolor}{rgb}{1.000000, 1.000000, 1.000000}
\pgfsetfillcolor{dialinecolor}
\pgfsetlinewidth{0.200000\du}
\pgfsetdash{{1.000000\du}{1.000000\du}}{0\du}
\pgfsetdash{{1.000000\du}{1.000000\du}}{0\du}
\pgfsetmiterjoin
\pgfsetbuttcap
\definecolor{dialinecolor}{rgb}{0.000000, 0.000000, 0.000000}
\pgfsetstrokecolor{dialinecolor}
\pgfpathmoveto{\pgfpoint{5.750000\du}{8.300000\du}}
\pgfpathcurveto{\pgfpoint{5.900000\du}{18.079203\du}}{\pgfpoint{36.200000\du}{16.729203\du}}{\pgfpoint{36.150000\du}{7.791703\du}}
\pgfpathcurveto{\pgfpoint{36.100000\du}{-1.145797\du}}{\pgfpoint{5.600000\du}{-1.479203\du}}{\pgfpoint{5.750000\du}{8.300000\du}}
\pgfusepath{stroke}
\pgfsetlinewidth{0.200000\du}
\pgfsetdash{}{0pt}
\pgfsetdash{}{0pt}
\pgfsetbuttcap
{
\definecolor{dialinecolor}{rgb}{0.000000, 0.000000, 0.000000}
\pgfsetfillcolor{dialinecolor}
\definecolor{dialinecolor}{rgb}{0.000000, 0.000000, 0.000000}
\pgfsetstrokecolor{dialinecolor}
\draw (17.845000\du,8.660000\du)--(13.345000\du,8.610000\du);
}
\pgfsetlinewidth{0.200000\du}
\pgfsetdash{}{0pt}
\pgfsetdash{}{0pt}
\pgfsetbuttcap
{
\definecolor{dialinecolor}{rgb}{0.000000, 0.000000, 0.000000}
\pgfsetfillcolor{dialinecolor}
\definecolor{dialinecolor}{rgb}{0.000000, 0.000000, 0.000000}
\pgfsetstrokecolor{dialinecolor}
\draw (20.745000\du,11.460000\du)--(20.750000\du,16.929203\du);
}
\definecolor{dialinecolor}{rgb}{1.000000, 1.000000, 1.000000}
\pgfsetfillcolor{dialinecolor}
\pgfpathellipse{\pgfpoint{20.745000\du}{8.610000\du}}{\pgfpoint{2.900000\du}{0\du}}{\pgfpoint{0\du}{2.800000\du}}
\pgfusepath{fill}
\pgfsetlinewidth{0.100000\du}
\pgfsetdash{}{0pt}
\pgfsetdash{}{0pt}
\definecolor{dialinecolor}{rgb}{0.000000, 0.000000, 0.000000}
\pgfsetstrokecolor{dialinecolor}
\pgfpathellipse{\pgfpoint{20.745000\du}{8.610000\du}}{\pgfpoint{2.900000\du}{0\du}}{\pgfpoint{0\du}{2.800000\du}}
\pgfusepath{stroke}
\definecolor{dialinecolor}{rgb}{0.000000, 0.000000, 0.000000}
\pgfsetstrokecolor{dialinecolor}
\node[anchor=west] at (18.5\du,9.\du){$\ten{A}$};
\definecolor{dialinecolor}{rgb}{1.000000, 1.000000, 1.000000}
\pgfsetfillcolor{dialinecolor}
\pgfpathellipse{\pgfpoint{10.445000\du}{8.610000\du}}{\pgfpoint{2.900000\du}{0\du}}{\pgfpoint{0\du}{2.800000\du}}
\pgfusepath{fill}
\pgfsetlinewidth{0.100000\du}
\pgfsetdash{}{0pt}
\pgfsetdash{}{0pt}
\definecolor{dialinecolor}{rgb}{0.000000, 0.000000, 0.000000}
\pgfsetstrokecolor{dialinecolor}
\pgfpathellipse{\pgfpoint{10.445000\du}{8.610000\du}}{\pgfpoint{2.900000\du}{0\du}}{\pgfpoint{0\du}{2.800000\du}}
\pgfusepath{stroke}
\definecolor{dialinecolor}{rgb}{0.000000, 0.000000, 0.000000}
\pgfsetstrokecolor{dialinecolor}
\node[anchor=west] at (8\du,9.\du){$\mat{B}$};
\pgfsetlinewidth{0.200000\du}
\pgfsetdash{}{0pt}
\pgfsetdash{}{0pt}
\pgfsetbuttcap
{
\definecolor{dialinecolor}{rgb}{0.000000, 0.000000, 0.000000}
\pgfsetfillcolor{dialinecolor}
\definecolor{dialinecolor}{rgb}{0.000000, 0.000000, 0.000000}
\pgfsetstrokecolor{dialinecolor}
\draw (7.545000\du,8.610000\du)--(2.945550\du,8.660550\du);
}
\pgfsetlinewidth{0.200000\du}
\pgfsetdash{}{0pt}
\pgfsetdash{}{0pt}
\pgfsetbuttcap
{
\definecolor{dialinecolor}{rgb}{0.000000, 0.000000, 0.000000}
\pgfsetfillcolor{dialinecolor}
\definecolor{dialinecolor}{rgb}{0.000000, 0.000000, 0.000000}
\pgfsetstrokecolor{dialinecolor}
\draw (28.645000\du,8.610000\du)--(23.645000\du,8.660000\du);
}
\definecolor{dialinecolor}{rgb}{1.000000, 1.000000, 1.000000}
\pgfsetfillcolor{dialinecolor}
\pgfpathellipse{\pgfpoint{31.545000\du}{8.610000\du}}{\pgfpoint{2.900000\du}{0\du}}{\pgfpoint{0\du}{2.800000\du}}
\pgfusepath{fill}
\pgfsetlinewidth{0.100000\du}
\pgfsetdash{}{0pt}
\pgfsetdash{}{0pt}
\definecolor{dialinecolor}{rgb}{0.000000, 0.000000, 0.000000}
\pgfsetstrokecolor{dialinecolor}
\pgfpathellipse{\pgfpoint{31.545000\du}{8.610000\du}}{\pgfpoint{2.900000\du}{0\du}}{\pgfpoint{0\du}{2.800000\du}}
\pgfusepath{stroke}
\definecolor{dialinecolor}{rgb}{0.000000, 0.000000, 0.000000}
\pgfsetstrokecolor{dialinecolor}
\node[anchor=west] at (30.\du,9.\du){$\mat{c}$};
\definecolor{dialinecolor}{rgb}{0.000000, 0.000000, 0.000000}
\pgfsetstrokecolor{dialinecolor}
\node[anchor=west] at (18.5\du,3.3\du){$\mat{D}$};
\end{tikzpicture}
\end{minipage}}
\caption{Basic TN diagrams.}
\label{fig:TNdiagrams}
\end{figure}

Expressions such as equation~\eqref{eq:TRdef} quickly become complicated as the number of tensors $d$ increases. Fortunately, there is a convenient graphical representation of these tensor structures. Figure~\ref{fig:TNdiagrams}(a) shows the diagram representation of a scalar, vector, matrix and 3-way tensor. Each of these tensors is represented by a node, and each edge corresponds with a particular index. Evidently, a scalar has no edges. The single most important operation with tensors is the summation over a particular index, also called a contraction, which is graphically represented by an edge that connects two nodes. For example, Figure~\ref{fig:TNdiagrams}(b) shows a diagram representation of a matrix $\mat{D}$ with entries
\begin{align}
    \mat{D}(i,j)&= \sum_{r_1=1}^{R_1} \sum_{r_2=1}^{R_2} \mat{B}(i,r_1) \; \ten{A}(r_1,j,r_2) \; \mat{c}(r_2), 
\label{eqn:simpleTN}
\end{align}
where the two contractions over the auxiliary indices $r_1,r_2$ are represented by the two edges connecting $\mat{B}$ with $\ten{A}$ and $\ten{A}$ with $\mat{c}$, respectively. 
Figure~\ref{def:TR} shows the diagram representation of a TR-decomposition of a $d$-way tensor $\ten{A}$ with uniform dimensions $N$. Fixing the indices of $\ten{A}$ to specific values removes all vertical edges in the diagram, which then represents equation~\eqref{eq:TRdefalt}. Assuming a uniform TR-rank $R$ sets the total storage complexity of the $d$ TR-tensors to $dNR^2$, completely removing the exponential dependence on $d$ and hence lifting the curse of dimensionality when $R \ll N^D$.
\begin{figure}[ht]
    \centering
    \input{figs/Curse.tex}
    \caption{The tensor $\ten{A}$ of size $N^d$ is represented by a TR $\ten{A}^{(1)},\ldots,\ten{A}^{(d)}$.}
    \label{fig:Curse}
\end{figure}
 
\section{Inability to recover the minimial-rank TR from a TR with sub-optimal ranks}
\label{sec:minrank}

In this section the main problem with the TR-decomposition is discussed through worked-out numerical experiments. We first provide the definition of a minimal-rank TR as given in~\cite[p.~5]{mickelin2018tensor}.
\begin{definition}
\label{def:minrank}
A vector $\mat{r} \triangleq \left(R_1,R_2,\ldots,R_d,R_1\right)$ is a minimal rank of a tensor $\ten{A}$ if (i) there exists a TR-representation of $\ten{A}$ with TR-ranks $\mat{r}$, and (ii) any other TR-rank $\mat{r}'$
of $\ten{A}$ satisfies $\mat{r} \leq \mat{r}'$ under
the elementwise inequality in $\mathbb{R}^{d+1}$.
\end{definition}
Note that it is perfectly possible that a minimal-rank TR turns out to be a TT. A first issue with recovering the minimal-rank TR occurs with the addition of TRs. Indeed, given a TR of $\ten{A}$ with TR-ranks $R_1,\ldots,R_d$, then the addition of $\ten{A}$ with itself results in a TR with ranks $2R_1,\ldots,2R_d$. Since the resulting tensor is merely a scaled version of the original tensor, one should be able to retrieve the original TR-ranks $R_1,\ldots,R_d$ through rounding. This is unfortunately not the case, an issue that was not discussed in~\cite{zhao2016tensor} and has only been recently resolved in~\cite{mickelin2018tensor} through an alternative definition of the addition operation. The proposed solution consists of treating the ``end-tensors" of the TR as TT-tensors when adding them. There is no specific reason why one would choose the ``end-tensors" of the TR. One can take any two consecutive TR-tensors in the ring and treat them as TT-cores when adding them together to retrieve the original TR-ranks with rounding. This solution is not satisfactory for two reasons: first, it breaks the ``symmetry" of the TR by treating two of its tensors differently when adding. Second, this ad-hoc solution provides no insight whether or how other important operations can be redefined such that TR-rounding retrieves the minimal-rank representation. Two such operations that are overlooked in~\cite{mickelin2018tensor} and remain problematic for rounding are contractions and Hadamard products. The motivation for each of these operations, together with simple examples where TR-rounding fails are discussed in the following two subsections.

\subsection{Matrix matrix multiplication in TR form}
\label{subsec:matrixmatrix}
\begin{figure}[tbh]
    \centering
    \input{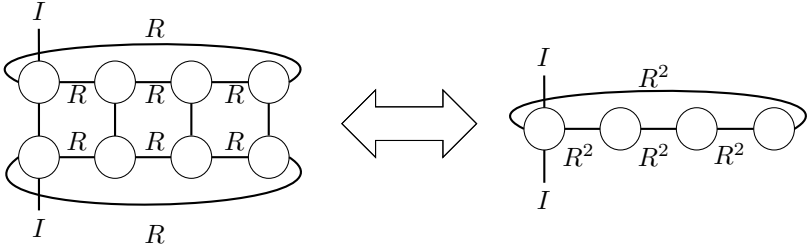}
    \caption{The matrix product $\mat{A}\mat{A}^T\in \mathbb{R}^{I \times I}$ in TR form when $d=4$.}
    \label{fig:mpocont}
\end{figure}

Repeated matrix matrix multiplications in TT/TR form are commonly used in e.g. the randomized subspace algorithm for computing a low-rank SVD approximation of a large matrix~\cite[p.~1235]{TNrSVD} or in the determination of the predicted covariance matrix in the tensor network Kalman filter~\cite{TNkalman,TNKalman2}. The modification of Definition~\ref{def:TR} for the TR/TT decomposition of a matrix $\mat{A} \in \mathbb{R}^{I_1I_2\cdots I_d \times J_1J_2\cdots J_d}$ is discussed in~\cite{Oseledets2010}. The main modification is that each TR-tensor $\ten{A}^{(k)} \in \mathbb{R}^{R_k \times I_k \times J_k \times R_{k+1}}$ is 4-way as it now contains both a row and column index. The entries of the $k$th TR-tensor that represents the matrix matrix multiplication $\mat{C}=\mat{A}\mat{B}$ for two matrices $\mat{A},\mat{B}$ represented by TR-tensors $\ten{A}^{(k)} \in \mathbb{R}^{R_k \times I_k \times J_k \times R_{k+1}},\,\ten{B}^{(k)} \in \mathbb{R}^{S_k \times J_k \times L_k \times S_{k+1}}$, respectively, are given by
\begin{align}
\label{eqn:matmatmul}
    \ten{C}^{(k)}([r_ks_k],i_k,l_k,[r_{k+1}s_{k+1}]) &= \sum_{j_k=1}^{J_k}\ten{A}^{(k)}(r_k,i_k,j_k,r_{k+1})\;\ten{B}^{(k)}(s_k,j_k,l_k,s_{k+1}).
\end{align}
The summation over all $k$ dimensions $J_k$ is effectively the summation over the columns of $\mat{A}$ and rows of $\mat{B}$ in the matrix matrix multiplication. The multi-indices $[r_ks_k]$, $[r_{k+1}s_{k+1}]$ can be converted to the linear indices through
\begin{align*}
    [r_ks_k] &= r_k + (s_k-1)\,R_k,\\
    [r_{k+1}s_{k+1}] &= r_{k+1} + (s_{k+1}-1)\,R_{k+1},
\end{align*}
which implies that the TR-ranks of $\ten{C}^{(k)}$ are the products of the corresponding TR-ranks of $\ten{A}^{(k)}$ with $\ten{B}^{(k)}$, which might be sub-optimal. A rounding step is therefore required in order to reduce the TR-ranks of $\mat{C}$ to their minimal values. The existence of a minimal-rank TR representation of $\mat{C}$ and the inability of retrieving this representation via TR-rounding can be demonstrated through the following example. Consider a rectangular matrix $\mat{A} \in \mathbb{R}^{I \times J^d }$ with $I \ll J^d$ that has an exact TR-representation with uniform TR-ranks $R$. The single row index of $\mat{A}$ is, without loss of generality, included into the TR-tensor $\ten{A}^{(1)} \in \mathbb{R}^{R_1 \times I \times J \times R_2}$ as an additional index. The matrix multiplication $\mat{A}\,\mat{A}^T \in \mathbb{R}^{I \times I}$ can now be computed according to~\eqref{eqn:matmatmul} in TR-form by summing over the $d$ indices with dimension $J$. These summations result in a ring of $d$ TR-tensors with uniform ranks $R^2$. The diagram of these contractions are shown in Figure~\ref{fig:mpocont} for the case $d=4$. The minimal-rank TR representation of $\mat{A}\,\mat{A}^T$ obtained from rounding is expected to be unit-rank, where the first TR-tensor is the $I \times I$ matrix $\mat{A}\,\mat{A}^T$ up to a scalar factor. The product of the remaining scalar TR-tensors constitute this factor.

\begin{figure}[ht]
    \centering
    \includegraphics[width=.85\textwidth]{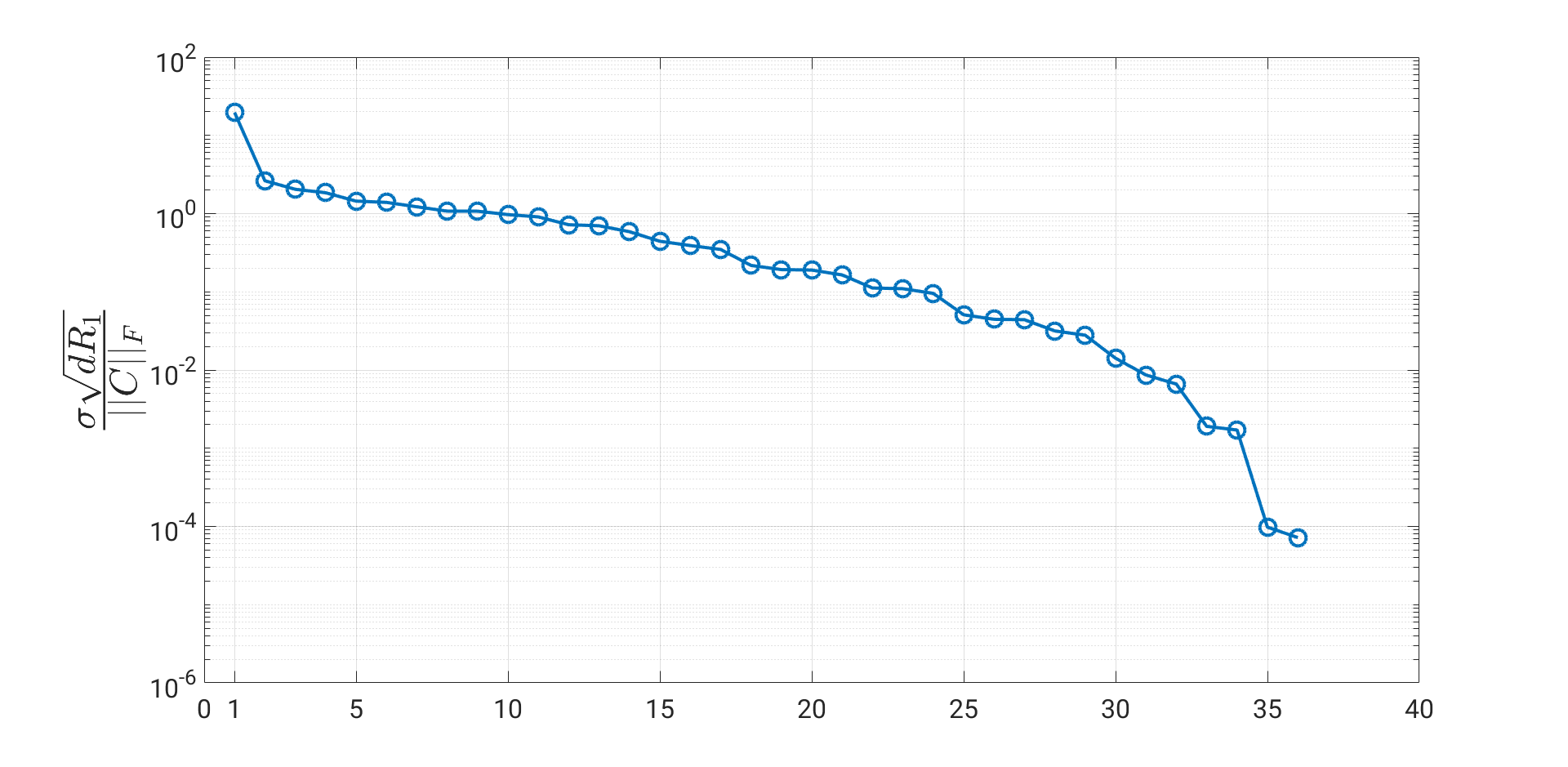}
    \caption{Scaled singular value profile of the unfolding of $\ten{C}^{(4)}$ when $R=6$.}
    \label{fig:singval}
\end{figure}

A numerical experiment is run on a TR that represents a $6 \times 6^4$ matrix $\mat{A}$ with uniform TR-ranks $R$ ranging from 1 up to 12. All entries of the TR-tensors are sampled from a standard normal distribution. Table~\ref{table:ATA} lists the retrieved ranks after rounding for varying $R$ values and shows that the TR-rounding algorithm (used with a tolerance of $10^{-10}$) fails to retrieve the minimal rank-1 representation. Indeed, no truncation of the TR-ranks is observed. The TR-rounding algorithm uses a tolerance $\delta = \epsilon ||\mat{C}||_F/\sqrt{dR_1}$ to determine the numerical rank of each unfolded TR-tensor via an SVD~\cite[p.~12]{mickelin2018tensor}. This tolerance $\delta$ guarantees that the relative error of the approximation is upper bounded by $\epsilon$. Figure~\ref{fig:singval} shows the singular value profile of the $\ten{C}^{(d)}$ tensor unfolded into a $36 \times 36$ matrix ($R=6$) in the TR-rounding algorithm, scaled by $\sqrt{dR_1}/||\mat{C}||_F$. From this figure one can deduce that the minimal relative truncation error $\epsilon$ is approximately $10^{-4}$. Furthermore, Figure~\ref{fig:singval} shows that a rank-1 truncation would result in a relative approximation error of 100\%, from which we can conclude that a TR simply does not allows us to retrieve the desired rank-1 representation.

The same matrix matrix multiplication can be computed using a TT representation by first converting the TR of $\mat{A}$ into a TT. The TT-rounding algorithm with a tolerance of $10^{-10}$ always retrieves the minimal rank-1 representation, as shown in the third column of Table~\ref{table:ATA}. The inability of the TR-rounding algorithm to reduce the obtained TR-ranks results in a linearly growing ratio of storage complexities, shown in the fourth column of Table~\ref{table:ATA} where the symbol $\#(\cdot)$ denotes the total number of parameters.
\begin{table}[tb]
\begin{center}
\caption{\label{table:ATA}Original TR-rank $R$ and rounded rank of $\mat{A}\mat{A}^T$ for both the TT and TR case in Figure~\ref{fig:mpocont}.}
\begin{tabular}{@{}rrrr@{}}
$R$  &  TR-round($R^2$) & TT-round($R^2$) & $\frac{\#(\textrm{TR}( \mat{A}\mat{A}^T))}{\#(\textrm{TT}( \mat{A}\mat{A}^T))}$ \\
          \midrule 
3	&   9 & 1 & 81 \\
6 	&  36 & 1 & 1296 \\
9   &  81 & 1 & 6561 \\
12  & 144 & 1 & 20736 
\end{tabular}
\end{center}
\end{table}

\subsection{Hadamard product of two tensors in TR form}
\label{subsec:hadamard}
The application of any nonlinear scalar function $f(\cdot)$ to all entries of a tensor $\ten{A}$ in TR/TT form can be approximated using Taylor series and implemented with Hadamard products. Suppose, without loss of generality, that the Mclaurin series can be used to approximate $f(\cdot)$, then
\begin{align*}
    f(\ten{A}) \approx f(0)\; \ten{I} + f'(0) \; \ten{A} + \frac{f''(0)}{2!} \; \left(\ten{A} \circ \ten{A} \right) + \cdots
\end{align*}
where $\ten{I}$ denotes the tensor of all ones of appropriate size. One application for such an approximation is the construction of a kernel matrix in support vector machines~\cite{SVDDV02} where the observed data matrix is given in TT/TR form~\cite{STTM}. Another use of the Hadamard product is found in the superfast Fourier Transform of data in TT form~\cite{dolgov2012superfast}. The TR that represents the Hadamard product $\ten{C}$ of two tensors $\ten{A},\, \ten{B}$ in TT form can be computed directly from the matrix Kronecker products
\begin{align}
\label{eq:hadamard}
\mat{C}(:,i_k,:) &= \mat{A}(:,i_k,:) \otimes \mat{B}(:,i_k,:) \; \; (1 \leq i_k \leq I_k, 1 \leq k \leq d),
\end{align}
as described in~\cite[p.~2310]{ivanTT}. Note that equation~\ref{eq:hadamard} was originally described for TTs but needs no modification to work for TRs. The Kronecker products in~\eqref{eq:hadamard} imply that the TR-ranks of $\ten{C}$ are the product of the TR-ranks of $\ten{A}$ with the corresponding TR-ranks of $\ten{B}$. After computing the TR of the Hadamard product, a rounding step can be used to reduce the ranks to their minimal values.

\begin{figure}[ht]
    \centering
    \includegraphics[width=.85\textwidth]{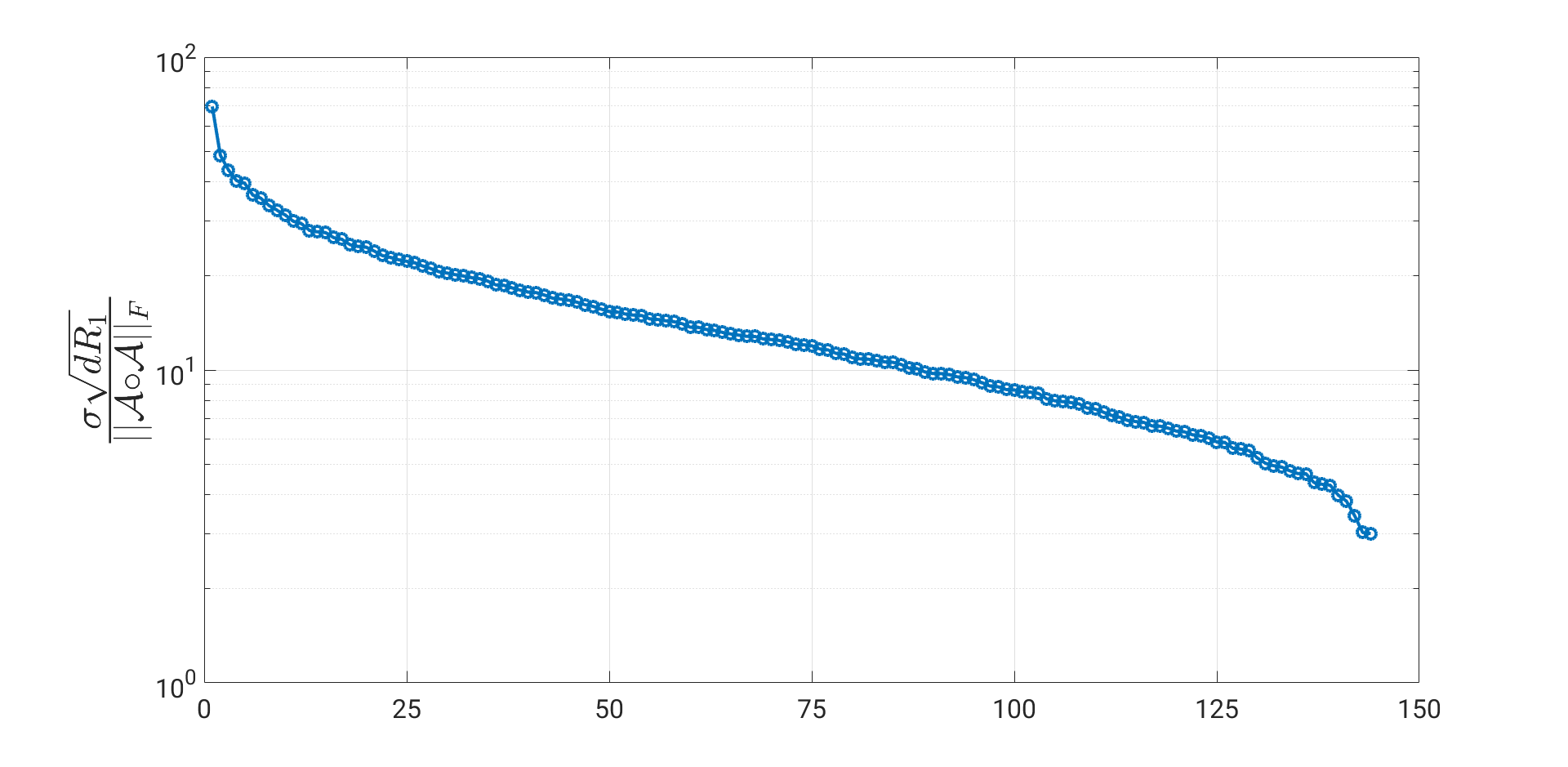}
    \caption{Scaled singular value profile of the unfolding of the fifth TR-tensor when $R=12$.}
    \label{fig:singval2}
\end{figure}

The inability of the TR-rounding algorithm to reduce the ranks of a Hadamard product is now demonstrated through a numerical experiment. A TR with uniform TR-ranks $R$ of a cubical 6-way tensor $\ten{A}$ with a dimension of 6 is hereto constructed. Again, each of the TR-tensor entries is sampled from a standard normal distribution. The Hadamard product of $\ten{A}$ with itself is then computed through~\eqref{eq:hadamard}, followed by the application of the TR-rounding algorithm. The whole procedure is repeated for TR-ranks ranging from 1 up to 12. Table~\ref{table:hadamard} lists the original TR-rank and the TR-rank obtained after rounding the Hadamard product. Just like in Section~\ref{subsec:matrixmatrix}, no reduction of the rank can be observed as each entry of the second column of Table~\ref{table:hadamard} is the square of the corresponding entry in the first column. Each of the TRs was also converted into a TT, which resulted in TT-ranks that were higher than the TR-ranks. The Hadamard product was computed with these TTs, followed by the TT-rounding algorithm. The maximal TT-ranks of the result are listed in the third column of Table~\ref{table:hadamard}. For $R\geq 5$ all maximal TT-ranks are 216, which implies that the total storage requirement of the Hadamard product in TT form remains constant for all $R\geq 5$. The inability of the rounding algorithm to reduce the storage complexity of the TR compared with the TT reflects itself in an increasing storage complexity, as shown in the fourth column of Table~\ref{table:hadamard}. One can see from the scaled singular value profile of the fifth TR-tensor of $\ten{A}\circ \ten{A}$ obtained through~\eqref{eq:hadamard} in Figure~\ref{fig:singval2} that no truncation of the rank is possible without the introduction of a relative approximation error of at least 100\%.  
\begin{table}[tb]
\begin{center}
\caption{\label{table:hadamard}Original TR-rank $R$ and rounded rank of the Hadamard product in both TT and TR form.}
\begin{tabular}{@{}rrrr@{}}
$R$& TR-round($R^2$) & max(TT-round($R^2$)) & $\frac{\#(\textrm{TR}(\ten{A}\circ \ten{A}))}{\#(\textrm{TT}( \ten{A}\circ \ten{A}))}$  \\
          \midrule 
3	& 9   &  45 & 0.13 \\
6 	& 36  & 216 & 0.48 \\
9   & 81  & 216 & 2.46\\
12  & 144 & 216 & 7.78
\end{tabular}
\end{center}
\end{table}
The existence of TR representations of $\ten{A}\circ \ten{A}$ with TR-ranks that are at least smaller than obtained through~\eqref{eq:hadamard}, but perhaps still not minimal, can be verified by converting the TT of $\ten{A}\circ \ten{A}$ into a TR and then applying the TR-rounding algorithm. In this case the rounding algorithm does reduce the TR-ranks and finds a maximal rank $R_5=108$ for all $R \geq 4$, such that lower-rank representations can be found for all $R \geq 11$. The TRs of $\ten{A}\circ \ten{A}$ obtained in this fashion exhibit smaller storage complexities than the TRs obtained through~\eqref{eq:hadamard} only when $R \geq 7$. 

\section{Conclusions}
\label{sec:conclusions}
The innate power and versatility of both TTs and TRs lies in their potential to perform linear algebra operations on vectors and matrices of exponential size in a computational and storage efficient manner. Numerical experiments in this article have demonstrated that it can be impossible to recover an exact minimal-rank TR representation through a rounding procedure after applying common operations such as contractions and the Hadamard product. As a result, TR-ranks will grow exponentially with repeated application of these operations. The impossibility of retrieving an exact minimal-rank TR through rounding, together with the non-existence of a best low-rank TR approximation severely limits the applicability of TRs.


\bibliographystyle{siamplain}
\bibliography{references}

\begin{thebibliography}{10}

\bibitem{TNkalman}
{\sc K.~Batselier, Z.~Chen, and N.~Wong}, {\em {A Tensor Network Kalman filter
  with an application in recursive MIMO Volterra system identification}},
  Automatica, 84 (2017), pp.~17--25.

\bibitem{TNKalman2}
{\sc K.~Batselier and N.~Wong}, {\em {Matrix output extension of the tensor
  network Kalman filter with an application in MIMO Volterra system
  identification}}, Automatica, 95 (2018), pp.~413--418.

\bibitem{TNrSVD}
{\sc K.~Batselier, W.~Yu, L.~Daniel, and N.~Wong}, {\em {Computing low-rank
  approximations of large-scale matrices with the Tensor Network randomized
  SVD}}, SIAM J. Matrix Anal. Appl., 39 (2018), pp.~1221--1244.

\bibitem{STTM}
{\sc C.~Chen, K.~Batselier, C.~Y. Ko, and N.~Wong}, {\em A support tensor train
  machine}, CoRR, abs/1804.06114 (2018), \url{http://arxiv.org/abs/1804.06114},
  \url{https://arxiv.org/abs/1804.06114}.

\bibitem{TNClassifier}
{\sc Z.~Chen, K.~Batselier, and N.~Wong}, {\em {Parallelized Tensor Train
  Learning of Polynomial Classifiers}}, IEEE Transactions on Neural Networks
  and Learning Systems, 29 (2018), pp.~4621--4632.

\bibitem{MAL-059}
{\sc A.~Cichocki, N.~Lee, I.~Oseledets, A.-H. Phan, Q.~Zhao, and D.~P. Mandic},
  {\em Tensor networks for dimensionality reduction and large-scale
  optimization: Part 1 low-rank tensor decompositions}, Foundations and
  Trends® in Machine Learning, 9 (2016), pp.~249--429.

\bibitem{MAL-067}
{\sc A.~Cichocki, A.-H. Phan, Q.~Zhao, N.~Lee, I.~Oseledets, M.~Sugiyama, and
  D.~P. Mandic}, {\em Tensor networks for dimensionality reduction and
  large-scale optimization: Part 2 applications and future perspectives},
  Foundations and Trends® in Machine Learning, 9 (2017), pp.~431--673.

\bibitem{dolgov2012superfast}
{\sc S.~Dolgov, B.~Khoromskij, and D.~Savostyanov}, {\em {Superfast Fourier
  transform using QTT approximation}}, Journal of Fourier Analysis and
  Applications, 18 (2012), pp.~915--953.

\bibitem{Espig2011}
{\sc M.~Espig, W.~Hackbusch, S.~Handschuh, and R.~Schneider}, {\em Optimization
  problems in contracted tensor networks}, Computing and Visualization in
  Science, 14 (2011), pp.~271--285.

\bibitem{Espig2012}
{\sc M.~Espig, K.~K. Naraparaju, and J.~Schneider}, {\em A note on tensor chain
  approximation}, Computing and Visualization in Science, 15 (2012),
  pp.~331--344.

\bibitem{handschuh2012changing}
{\sc S.~Handschuh}, {\em Changing the topology of tensor networks}, arXiv
  preprint arXiv:1203.1503,  (2012).

\bibitem{Khoromskij2011}
{\sc B.~N. Khoromskij}, {\em O(dlog{\thinspace}n)-quantics approximation of n-d
  tensors in high-dimensional numerical modeling}, Constructive Approximation,
  34 (2011), pp.~257--280.

\bibitem{Landsburg2012}
{\sc J.~M. Landsburg, Y.~Qi, and K.~Ye}, {\em On the geometry of tensor network
  states}, Quantum Info. Comput., 12 (2012), pp.~346--354.

\bibitem{LGROL14}
{\sc V.~Lebedev, Y.~Ganin, M.~Rakhuba, I.~Oseledets, and V.~Lempitsky}, {\em
  Speeding-up convolutional neural networks using fine-tuned cp-decomposition},
  arXiv preprint arXiv:1412.6553,  (2014).

\bibitem{Lee2015}
{\sc N.~Lee and A.~Cichocki}, {\em {Estimating a Few Extreme Singular Values
  and Vectors for Large-Scale Matrices in Tensor Train Format}}, SIAM J. Matrix
  Anal. Appl., 36 (2015), pp.~994--1014.

\bibitem{Lee2016}
{\sc N.~Lee and A.~Cichocki}, {\em {Regularized Computation of Approximate
  Pseudoinverse of Large Matrices Using Low-Rank Tensor Train Decompositions}},
  SIAM J. Matrix Anal. Appl., 37 (2016), pp.~598--623.

\bibitem{mickelin2018tensor}
{\sc O.~Mickelin and S.~Karaman}, {\em Tensor ring decomposition}, arXiv
  preprint arXiv:1807.02513,  (2018).

\bibitem{Novikov2015}
{\sc A.~Novikov, D.~Podoprikhin, A.~Osokin, and D.~Vetrov}, {\em {Tensorizing
  Neural Networks}}, in Proceedings of the 28th International Conference on
  Neural Information Processing Systems, NIPS'15, Cambridge, MA, USA, 2015, MIT
  Press, pp.~442--450.

\bibitem{Oseledets2010}
{\sc I.~Oseledets}, {\em {Approximation of $2^d \times 2^d$ Matrices Using
  Tensor Decomposition}}, SIAM J. Matrix Anal. Appl., 31 (2010),
  pp.~2130--2145.

\bibitem{Os-mvk2-2011}
{\sc I.~Oseledets}, {\em {DMRG} approach to fast linear algebra in the
  {TT}--format}, Comput. Meth. Appl. Math., 11 (2011), pp.~382--393,
  \url{https://doi.org/10.2478/cmam-2011-0021}.

\bibitem{ivanTT}
{\sc I.~Oseledets}, {\em {Tensor-Train Decomposition}}, SIAM J. Sci. Comput.,
  33 (2011), pp.~2295--2317, \url{https://doi.org/10.1137/090752286}.

\bibitem{Oseledets2012}
{\sc I.~Oseledets and S.~Dolgov}, {\em Solution of linear systems and matrix
  inversion in the {TT}-format}, SIAM Journal on Scientific Computing, 34
  (2012), pp.~A2718--A2739.

\bibitem{ttcross}
{\sc I.~Oseledets and E.~Tyrtyshnikov}, {\em {TT}-cross approximation for
  multidimensional arrays}, Linear Algebra and its Applications, 422 (2010),
  pp.~70--88.

\bibitem{MPS1997}
{\sc S.~Rommer and S.~\"Ostlund}, {\em Class of ansatz wave functions for
  one-dimensional spin systems and their relation to the density matrix
  renormalization group}, Phys. Rev. B, 55 (1997), pp.~2164--2181.

\bibitem{SCHOLLWOCK201196}
{\sc U.~Schollw\"ock}, {\em The density-matrix renormalization group in the age
  of matrix product states}, Annals of Physics, 326 (2011), pp.~96 -- 192.
\newblock January 2011 Special Issue.

\bibitem{stoudenmire2018learning}
{\sc E.~Stoudenmire}, {\em Learning relevant features of data with multi-scale
  tensor networks}, Quantum Science and Technology, 3 (2018), p.~034003.

\bibitem{stoudenmire2016supervised}
{\sc E.~Stoudenmire and D.~Schwab}, {\em Supervised learning with tensor
  networks}, in Advances in Neural Information Processing Systems, 2016,
  pp.~4799--4807.

\bibitem{SVDDV02}
{\sc J.~A.~K. Suykens, T.~{Van Gestel}, J.~{De Brabanter}, B.~{De Moor}, and
  J.~Vandewalle}, {\em Least Squares Support Vector Machines}, World
  Scientific, Singapore, 2002.

\bibitem{Wang_TRcompletion}
{\sc W.~Wang, V.~Aggarwal, and S.~Aeron}, {\em Efficient low rank tensor ring
  completion}, in 2017 IEEE International Conference on Computer Vision (ICCV),
  Oct 2017, pp.~5698--5706.

\bibitem{wang2017efficient}
{\sc W.~Wang, V.~Aggarwal, and S.~Aeron}, {\em Efficient low rank tensor ring
  completion}, in Proceedings of the IEEE Conference on Computer Vision and
  Pattern Recognition, 2017, pp.~5697--5705.

\bibitem{ye2018tensor}
{\sc K.~Ye and L.-H. Lim}, {\em Tensor network ranks}, arXiv preprint
  arXiv:1801.02662,  (2018).

\bibitem{yuan2018rank}
{\sc L.~Yuan, C.~Li, D.~Mandic, J.~Cao, and Q.~Zhao}, {\em Rank minimization on
  tensor ring: A new paradigm in scalable tensor decomposition and completion},
  arXiv preprint arXiv:1805.08468,  (2018).

\bibitem{zhao2016tensor}
{\sc Q.~Zhao, G.~Zhou, S.~Xie, L.~Zhang, and A.~Cichocki}, {\em Tensor ring
  decomposition}, arXiv preprint arXiv:1606.05535,  (2016).

\end{thebibliography}
\end{document}